\documentclass[preprint,12pt]{elsarticle}

\usepackage{amsmath,amssymb,mathdots,mathrsfs,mathtools}
\usepackage{bm,upgreek,nicefrac}
\usepackage[mathscr]{eucal}
\usepackage{xcolor}
\usepackage{tikz,accents,subcaption}

\newcommand{\bi}[1]{\emph{\textbf{#1}}}

\newcommand{\undb}[1]{\underaccent{\bar}{#1}}

\newtheorem{thm}{Theorem}[section]
\newtheorem{cor}[thm]{Corollary}
\newtheorem{lem}[thm]{Lemma}

\newtheorem{exmp}{Example}[section]
\newproof{pf}{Proof}

\journal{IET}
\date{}
\begin{document}
\begin{frontmatter}

\title{Generalizing Laplacian Controllability of Paths}

\author{Shun-Pin Hsu\corref{cor1}}
\cortext[cor1]{corresponding author}
\ead{shsu@nchu.edu.tw}\author[]{Ping-Yen Yang}
\ead{ d105064202@mail.nchu.edu.tw}


\address{Department of Electrical Engineering, National Chung Hsing University\\
250, Kuo-Kuang Rd., Taichung 402, Taiwan}

\begin{abstract}
It is well known that if a network topology is a path or line and the states of vertices or nodes evolve according to the consensus policy, then the network is Laplacian controllable by an input connected to its terminal vertex. In this work a path is regarded as the resulting graph after interconnecting a finite number of two-vertex antiregular graphs and then possibly connecting one more vertex. It is shown that the single-input Laplacian controllability of a path can be extended to the case of interconnecting a finite number of $k$-vertex antiregular graphs with or without one more vertex appended, for any positive integer $k$. The methods to interconnect these antiregular graphs and to select the vertex for connecting the single input that renders the network Laplacian controllable are presented as well.
\end{abstract}

\end{frontmatter}

\section{Introduction}
\begin{figure*}[t!]
		\centering
		\includegraphics[width=\linewidth]{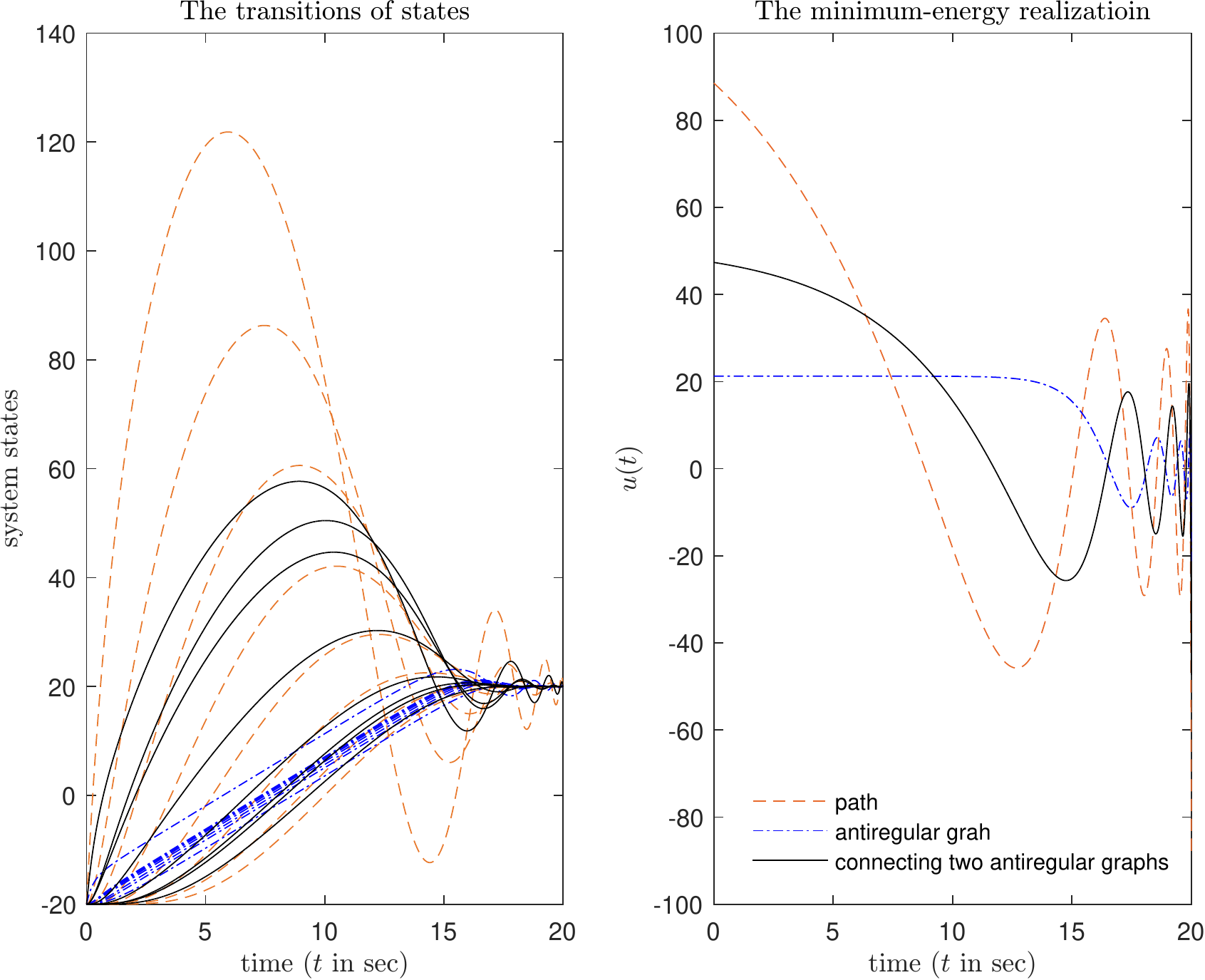}
		\caption{The transition and input comparisons of driving three different eight-state systems in (\ref{eq:mcons}) from $-20[\,1\,1\,1\,1\,1\,1\,1\,1\,\,]$ to $20[\,1\,1\,1\,1\,1\,1\,1\,1\,\,]$ under the minimum energy control~\cite[p.149]{chen99}. In the left subplot, the transitions of the eight system states defined on an eight-vertex path are shown in eight dashed lines. Those on an eight-vertex antiregular graph and on a graph constructed by interconnecting two four-vertex antiregular graphs are in eight dash-dot lines and in eight solid lines respectively. The right subplot shows the corresponding input signals that realize the minimum energy control. The comparisons show that the transition smoothness and the energy consumption in the case of interconnecting two four-vertex antiregular graphs, are between that of a eight-vertex path and that of a eight-vertex antiregular graph (cf. Fig.1 in~\cite{hsu17b}).}\label{fig:comp}
\end{figure*}

Integrating modern sensing, communication, and control technologies into conventional systems to create their state-of-the-art versions is the recent trend in industry. These innovative systems are able to collect real-time data, analyze the environment parameters or background signals, and respond adaptively. Immediate benefits from these modern systems might include the high production efficiency, low maintenance costs and so on. An interesting example is called the energy internet. This technology is based on the integration of power flow and information flow, and serves as a response to the call for effective and efficient operations of the power grids established by myriads of homes, retail stores and factories. Clearly, as the power generation mode is gradually transformed from the centralized mass power plants to distributed micro power plants, various issues such as the coordination control, network security, smart metering and operation management need to be addressed ~\cite{pat16,sun16,kal16,mon16,lin17,han18,lin18}. In fact, many of these issues are closely related to the mechanism of information dissemination on the network, and have their impact on the overall performance. For instance, in large solar power or wind farms, solar panels or wind turbines might need to track the elevation angle of the sun or incident angle of the airflow to maximize the efficiency of electricity generation. If each power generation unit is controlled by independent signals for the angle tracking, the activation cost will be very high as the number of units is large. It is expected that the control signals can be sent through some circuit that is appropriately designed so that only a few control inputs are needed. This circuit propagates the control signals provided by several input sources and drives all distributed power generation units to reach their respective states. Unfortunately, even for a relative simple circuit that can be described by a linear and time-invariant system, finding a small set of variables to affect with an input to achieve the so-called minimal controllability is a NP-hard problem~\cite{ols14}.
\begin{figure*}[t!]
		\centering
		\includegraphics[width=\linewidth]{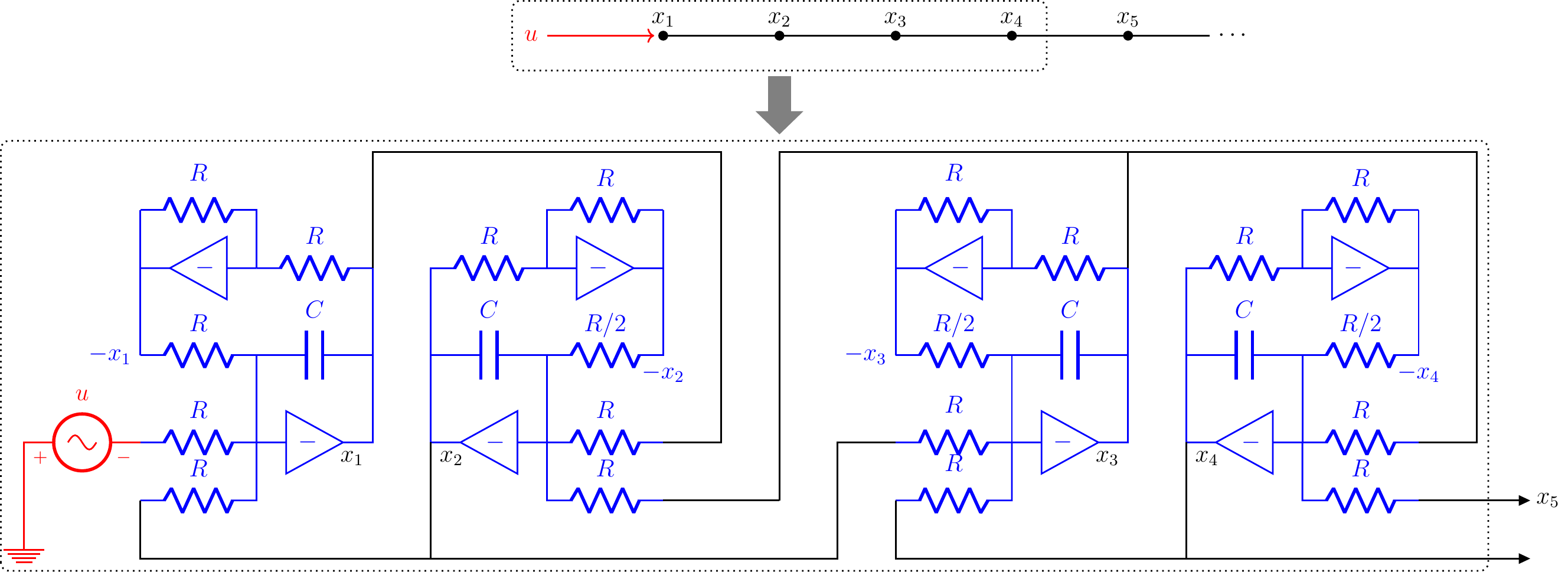}
		\caption{A circuit implementation of a Laplacian controllable path, where $RC=1$.}\label{fig:path}
	\end{figure*}
To simplify the broadcasting mechanism of control signals, we expect the control information to be exchanged only locally but gradually spread all over the entire system. The advantages of this mechanism include the simplicity of the circuit and the state consistency of the zero-input response. Moreover, the problem can be formulated in the context of Laplacian controllability of a multiagent system with leader-follower Laplacian dynamics~\cite{agu15}. An example of the circuit is shown in~Fig.\ref{fig:path}, where only one input is used. The structure is composed of nearly cascaded connection of blocks of electric elements and follows the network topology of a path (or called a line) graph. A simple spectra analysis of this path reveals that it is Laplacian controllable by an input connected to its terminal vertex~\cite[Lemma 3.2]{par12}. That is, starting from any state $\bi{x}=[\,x_1\,x_2,\cdots\,]^T$, a control function $\bi{u}$ capable of driving $\bi{x}$ to any specified $\bi{x}$' at any specified finite time exists. The Laplacian controllability problem of linear systems has been an active research topic in the past decade due to its important application in multiagent networks~\cite{mes10}. A standard approach to attacking the problem is to analyze the properties of Laplacian eigenvectors of the graphs and apply the Popov-Belevitch-Hautus (PBH) test to draw the conclusion~\cite{rah09}. It is easier to derive the conditions that lead to uncontrollability than to controllability since an apparent class of uncontrollable graphs can be identified from the symmetry of network topology, in some sense~\cite{eg12}. For a general graph with hundreds of vertices, a simple method to decide the minimum number of inputs and the vertices to connect these inputs to ensure the Laplacian controllability is yet unknown. Although simple algorithms based on vertex partitioning to output a lower bound and an upper bound for the minimum number of controllers have been proposed~\cite{zha14}, the general gap between the bounds needs future improved schemes to become small. If the network topology follows specific patterns, the Laplacian eigespaces might become tractable, and thus determining the minimum number of controllers and locating the vertices to connect these controllers might be possible~\cite{not13,nab13,agu15b,hsu17a}. Among these results, only a few types of graphs can be controlled by one input. Recently, two new members in the class of single-input Laplacian controllable graphs were discovered. An important insight into the discovery is that interconnecting two single-input Laplacian controllable graphs appropriately might preserve the single-input controllability. The first member features the interconnection of a path and an antiregular graph~\cite{hsu17b}, and the second, the interconnection of two antiregular graphs~\cite{hsu18}. An antiregular graph is a connected simple graph that has exactly one pair of degree-repeating vertices, or the vertices that have the same number of neighboring vertices~\cite{mer03}. In Fig.~\ref{fig:comp} we compare the transient responses and the minimum energy required to drive three eight-vertex single-input Laplacian controllable graphs. The figure illustrates that a path is more difficult to drive than an antiregular graph and the transient behavior of states with the connection topology of a path is less smooth. However, an antiregular graph requires a dominating vertex to connect all other vertices. This requirement becomes infeasible or impractical as the graph has a large number of vertices. The third graph is a constructed by interconnecting two four-vertex antiregular graphs. Unsurprisingly, the figure shows that its performance is somewhere between the first two cases. In this work we follow the line of our previous study to explore more members in the class of single-input Laplacian controllability. Firstly, we take a closer look at the simple structure of a $v$-vertex path and interpret its connecting style as the interconnection of a finite number of small-sized antiregular graphs. Specifically, if $v$ is odd, we start from an isolated vertex and connect this vertex to the degree-repeating vertex of a two-vertex antiregular graph. If $v$ is even, we simply start from a two-vertex antiregular graph. Note that at this moment the graph is Laplacian controllable by an input connected to the starting vertex as $v$ is odd, or by an input connected to the degree-repeating vertex as $v$ is even. Each time a new two-vertex antiregular graph is added, the degree-repeating vertex of the new graph is connected to the terminal vertex of existing interconnected graph. Clearly, the resulting graph after the addition is still single-input Laplacian controllable by the input connected to the vertex aforementioned. With this interpretation, we show if we increase the number of vertices of the antiregular graphs from two to any finite integer, the resulting graph is still Laplacian controllable by the input specified above.
The major contribution of our result is the generalization of the single-input Laplacian controllability of a path. An immediate benefit of this generalization is the assurance of one more member in the class of single-input Laplacian controllable graphs. The existence of this member is not verified through the full knowledge of Laplacian eigespaces of the graphs, but through the indirect approach that proves the distinctness of every Laplacian eigenvalue. This approach demonstrates the possibility to identify a single-input Laplacian controllable graph whose Laplacian eigenspaces are only partially tractable. The second benefit of this generalization is to offer more options for the design of network topology subject to some edge constraints or control costs. As illustrated in Fig.~\ref{fig:comp}, different connection topologies take different amounts of energy to drive and have different levels of transient fluctuations. These differences are closely related to the edge parameters such as the diameter and the maximum number of edges connected to one vertex. Our result allows one to construct two types of single-input Laplacian controllable graphs. The first type has vertex number $kn$, where $k\geq 3,n\geq 1$, diameter $3n-1$, and the maximum number of edges connected to one vertex $k-1$. The corresponding values of the second type are $kn+1$, $3n$, and $k-1$ respectively. Observe in the case with $kn$ vertices, a path and an antiregular graph have diameters $kn-1$ and $2$ respectively, and the maximum number of edges connected to one vertex $2$ and $kn-1$ respectively. We can apply our result to design a Laplacian controllable and feasible connection topology by adjusting $k$ and $n$. This flexibility is not available using a path or an antiregular graph alone.

The rest of this paper is organized as follows. The notations and essential concepts used throughout the paper are reviewed in the second section. Our main results on generalizing the Laplacian controllability of path graphs are presented in the third section. For better presentation, the results are accompanied with numerical examples. The paper is concluded in Section~IV where some future research topics of interest are discussed.

\section{Preliminaries}
We first define some symbols and recapitulate important concepts to be used in deriving our main results. More details can be found in the standard textbooks dealing with similar topics or articles referred. Let $\mathbb{R}^d$ be the set of $d$-entry real column vectors and $\mathbb{R}^{d_1\times d_2}$ the set of real matrices with size $d_1\times d_2$. A matrix $A\in\mathbb{R}^{d\times d}$ is called a square matrix of order $d$. We use $\textbf{1}$ and $\textbf{0}$ to represent the column vectors of 1's and 0's respectively and $I_k$ the identity matrix of order $k$. The $i$th column of an identity matrix is written as $\bi{e}_i$. Its size is subject to the context. For two sets $S_1$ and $S_2$, the set difference $S_1\setminus S_2$ is defined as $\{s|s\in S_1, s\notin S_2\}$. We use $\lfloor x\rfloor$ and $\lceil x\rceil$ to represent the largest integer not greater than $x$ and the smallest integer not less than $x$, respectively. We say $(\lambda,\bi{v})$ is an eigenpair of $P$ if $\bi{v}$ is a corresponding eigenvector of $P$ to the eigenvalue $\lambda$. The Kronecker product of matrices $P_1$ and $P_2$ is written as $P_1\otimes P_2$. A $k$-vertex graph can be defined by a two-tuple $(V,E)$ where $V:=\{1,2,\cdots,k\}$ is called the vertex set and $E:=\{\,(v_1,v_2)\,|\,v_1,v_2\in V, v_1\neq v_2\}$ the edge set. A graph is called a simple graph if it is undirected and unweighted. Thus for a simple graph, any entry in the edge set is an unordered pair, namely, $E:=\{\,\{v_1,v_2\}\,|\,v_1,v_2\in V, v_1\neq v_2\}$, and all edges defined in $E$ have the same weight. In this work we restrict our discussion to connected simple graphs only. More details on the algebraic aspects of such graphs can be seen, for example, in~\cite{god01,biggs93}. An example of connected simple graph is the path. A $k$-vertex path has the vertex set $V=\{1,2,\cdots,k\}$ and edge set $E=\{\,\{1,2\},\{2,3\},\cdots,\{k-1,k\}\}$. Two vertices $v_1$ and $v_2$ in $V$ are neighbors if $\{v_1,v_2\}\in E$. The neighbor set of vertex $v$ is $\mathcal{N}_v:=\{u\,|\,\{v,u\}\in E\}$. The number of elements in $\mathcal{N}_v$, written as $|\mathcal{N}_v|$, is called the degree or valency of vertex $v$. A connected simple graph has at least two vertices sharing the same degree. Those vertices are called degree-repeating vertices. A vertex is called a terminal vertex if it has only one neighbor. A vertex is called a dominating vertex or universal vertex if all other vertices in the graph are its neighbors. We can number the vertices of a $k$-vertex graph such that $d_i$, the degree of the $i$th vertex, satisfies $d_i\geq d_{i+1}$ for each $i\in\{1,2,\cdots,k-1\}$. The sequence $\bi{d}:=(\,d_1,d_2,\cdots,d_k\,)$ is called the degree sequence of the graph and $\tau_{\bi{d}}$, the trace of $\bi{d}$, is defined as $|\,j:d_j\geq j\,|$. The conjugate of $\bi{d}$ is $\bi{d}^*:=\left(\,d_1^*,d_2^*,\cdots,d_k^*\,\right)$  where $d_i^*=|\,j:d_j\geq i\,|$. A sequence $\bi{a}=(\,a_1,a_2,\cdots,a_k\, )$ is called graphical if there exists a $k$-vertex graph whose degree sequence is exactly $\bi{a}$. It was shown~\cite[p.72]{mol12} that the necessary and sufficient condition for $\bi{d}$ to be graphical is that
\begin{equation}\label{lem:gra}
\sum_{i=1}^j(d_i+1)\leq\sum_{i=1}^jd_i^*, \quad\mbox{$\forall j\in\{1,2,\cdots,\tau_\bi{d}\}$}.
\end{equation}
Suppose $(V,E)$ determines a connected simple graph and its degree sequence is $\bi{d}=(\,d_1,d_2,\cdots,d_k\, )$. The Laplacian matrix $\mathcal{L}$ of the graph is defined as
\begin{equation*}
\mathcal{L}:=\mathcal{D}-\mathcal{A}
\end{equation*}
where $\mathcal{D}$ is a diagonal matrix whose $i$th diagonal term is $d_i$, and $\mathcal{A}$ is a binary matrix whose $(i,j)$th element is $1$ if $\{i,j\}\in E$ and is $0$ otherwise.
As an example, the Laplacian matrix $\mathcal{L}$ of a $k$-vertex path graph can be written as
\begin{equation*}
\mathcal{L}=\begin{bmatrix*}[r]
1 &  - 1   &&&&     \\
-1   & 2 &-1 &&&        \\
 &-1  &2 &-1   & &        \\
  &&\ddots &\ddots &\ddots &\\
 &&&-1 &2 &-1\\
  & &&&-1&1
\end{bmatrix*}_{k\times k}.
\end{equation*}
The eigenvalues and eigenvectors of $\mathcal{L}$ are called the Laplacian eigenvalues and Laplacian eigenvectors, respectively, of the graph determined by $(V,E)$. Clearly $(0,\textbf{1})$ is an eigenpair of $\mathcal{L}$, and all eigenvectors of $\mathcal{L}$ are orthogonal to $\textbf{1}$ except for those corresponding to the zero eigenvalue. More properties concerning the eigenvalues and eigenvectors of $\mathcal{L}$ were summarized in~\cite{mer98}.
 According to the Grone-Merris theorem, the spectrum of $\mathcal{L}$ is majorized by the conjugate $\bi{d}^*$ of the degree sequence of the graph, namely,
\begin{equation}\label{lem:GroMer}
\sum_{i=1}^t\ell_{k-i+1}\le\sum_{i=1}^td_i^*\quad\forall t\in\{1,2,\cdots,k\}
\end{equation}
where $\ell_i$ is the $i$th smallest Laplacian eigenvalue of the graph. In the special case that the equality in~(\ref{lem:gra}) holds, $\bi{d}$ determines uniquely a \emph{threshold graph} or called \emph{maximal graph}, which turns out to be the special case that the equality in~(\ref{lem:GroMer}) holds~\cite{mer94}. Namely, the Laplacian eigenvalues of threshold graphs are readily available from their degree sequences. In fact, threshold graphs admit several equivalent definitions to the version aforementioned. For example, they can be defined via the symmetry of the Ferrers-Sylvester diagram~\cite[p.70]{mol12}, or defined through the graph constructing process involving join and union operations only~\cite{bap13}. An antiregular graph is a connected simple graph that has exactly one pair of vertices sharing the same degree~\cite{mer03}. It can be easily verified that an antiregular graph belongs to the class of threshold graphs. Elegant Laplacian eigenpair properties of threshold graphs can thus be applied directly.

An autonomous system that is linear, time-invariant and evolves according to the consensus policy
has the following form:
\begin{equation}\label{eq:cons}
\dot{x}_i=-\sum_{j\in\mathcal{N}_i}(x_i-x_j),
\end{equation}
where $\mathcal{N}_i$ is a subset of the set of state variables. This system is usually employed to model a dynamic network system in which each vertex state interacts with its neighboring vertex states for communication such that the local information can be propagated throughout the entire system. Adopting the Laplacian matrix of the graph that models the network system to express~(\ref{eq:cons}) yields the Laplacian dynamics~\cite[p.1613]{agu15}
\begin{equation}\label{eq:automcons}
\dot{\bi{x}}=-\mathcal{L}\bi{x}.
\end{equation}
To control the autonomous system in~(\ref{eq:automcons}) with one controller, we can apply input $u(t)$ via a binary control vector $\bi{b}\in\{0,1\}^{k}$ so that
\begin{equation}\label{eq:mcons}
\dot{\bi{x}}=-\mathcal{L}\bi{x}+\bi{b}u(t)
\end{equation}
where the $i$th element of $\bi{b}$ is
$1$ if vertex $i$ is connected to input $u(t)$, and is $0$ otherwise. For simplicity, we use the notation $(\mathcal{L},\bi{b})$ to represent the controlled graph model of the dynamic system in~(\ref{eq:mcons}). We say a graph is single-input Laplacian controllable by the input $u$ if its corresponding $(\mathcal{L},\bi{b})$ is controllable. A well-known result on a path graph is that it is Laplacian controllable by an input connected to one of its terminal vertex. In this paper we first interpret a $v$-vertex path, where $v$ is even, as an interconnection of a finite number of two-vertex antiregular graphs. In the case that $v$ is odd, one more vertex is needed for the interconnection. We show that the Laplacian controllability is preserved if we increase the number of vertices of the antiregular graph from two to any finite number and the vertices for interconnecting the antiregular graphs or for connecting the input are appropriately selected. Our result is based on the following version of Popov-Belevitch-Hautus theorem for symmetric matrices. It is a standard tool for the analysis of Laplacian controllability of a graph.
 \begin{thm}\label{thm:ns}~(cf. \cite[p.145]{chen99}) A graph is Laplacian controllable by an input via a control vector,  if and only if the graph does not have a Laplacian eigenvector orthogonal to the control vector.
 \end{thm}

\section{Main Results}
Let $\mathcal{C}_A$ be the class of antiregular graphs defined above. A $k$-vertex graph in $\mathcal{C}_A$ is written as $\mathbb{G}^{(k)}_A$. The Laplacian matrix corresponding to $\mathbb{G}^{(k)}_A$ is denoted by $\mathcal{L}^{(k)}_A$. Let the $i$th entry in the degree sequence be the $i$th diagonal term. We obtain the $\mathcal{L}^{(k)}_A$ written as
\[
\begin{bmatrix}
k-1 &    -1   & \cdots &    -1   &    -1   & \cdots &    -1   &    -1   \\
-1   & k-2 & \cdots &    -1   &    -1   & \cdots &    -1   &        \\
\vdots & \vdots & \ddots & \vdots & \vdots &\iddots        &        &        \\
-1   &    -1   & \cdots & -\underline{\kappa} & \beta_k & & & \\
-1   &    -1   & \cdots & \beta_k  & -\underline{\kappa} & & & \\
\vdots & \vdots &\iddots        &        &        & \ddots &        &        \\
-1   &    -1   &        &        &        &        &   2   &        \\
-1   &        &        &        &        &        &        &   1
\end{bmatrix}
\]
where $\beta_k$ is $1$ if $k$ is even and is $0$ if $k$ is odd.
\begin{equation}\label{def:undbk}
\undb{\kappa}:=\left\lfloor\frac{k}{2}\right\rfloor.
\end{equation}
Note that an antiregular graph is a special threshold graph, which can be constructed by adding vertices one by one
starting from an isolated vertex. Each time a new vertex is added, only the union or join operation can be applied. This constructing process leads to an integral Laplacian spectrum that can be obtained from the degree sequence, and a set of orthogonal Laplacian eigenvectors that can be derivable from the Laplacian matrix in a straightforward method~\cite{mer98}. In case of a $k$-vertex antiregular graph, the set of its Laplacian eigenvalues can be written as
\begin{equation}\label{def:Lda}
	\Lambda_k:=\{0,1,\cdots,k\}\setminus\{\bar{\kappa}\}.
\end{equation}
where
\begin{equation}\label{def:bk}
\bar{\kappa}=\left\lceil\frac{k}{2}\right\rceil,
\end{equation}
and a full set of orthogonal eigenvectors of $\mathcal{L}^{(k)}_A$ can be obtained using the following lemma.
\begin{lem}(\cite{hsu16,agu15b}) \label{lem:ev}
	Let the $(i,j)$th entry of matrices $T^{(m)}$ be $t_{ij}^{(m)}$ for each $m\in\{1,2,3,4\}$. Suppose $T^{(1)}=\mathcal{L}^{(k)}_A$ and let $T^{(2)},T^{(3)}$ be generated by
	\begin{equation} \label{def:t2}
	t_{ij}^{(2)}=\left\{
	\begin{array}{cl}
	-1-t_{ij}^{(1)}, & \mbox{if $ j>i$} \\
	t_{ij}^{(1)},    & \mbox{o.w.}
	\end{array}
	\right.
	\end{equation}
	and
	\begin{equation} \label{def:t3}
	t_{ij}^{(3)}=\left\{
	\begin{array}{cl}
	-\sum_{k,k\ne j} t_{kj}^{(2)}, & \mbox{if  $j=i$} \\
	t_{ij}^{(2)},                  & \mbox{o.w.}
	\end{array}
	\right..
	\end{equation}
Finally, remove the (unique) zero column of $T^{(3)}$ and append the column of $1$'s (or $-1$'s) to the last column to yield $T^{(4)}$. Then the $j$th column of $T^{(4)}$ is the eigenvector corresponding to the $j$th largest eigenvalue of $\mathcal{L}^{(k)}_A$, or the $j$th entry of the conjugate of the degree sequence of $\mathbb{G}^{(k)}_A$.
\end{lem}
Now we consider the interconnection of $n$ antiregular graphs $\mathbb{G}_A^{(k)}$.
The interconnection is initiated with a $k$-vertex antiregular graph where $k\ge 2$. To add  a new antiregular graph, a new edge is employed to connect one of the two vertices with the same degree in the new antiregular graph, and the terminal vertex of the existing graph. The resulting graph after interconnecting $n$ $k$-vertex antiregular graphs has the Laplacian matrix
\begin{align*}
\mathscr{L}_n^{k}:&=\mathcal{L}^k_n+ZZ^T\\
&=\mathcal{L}^k_n+\sum_{i=1}^{n-1}\bi{z}_i\bi{z}_i^T
 \end{align*}
where
\begin{align*}
\mathcal{L}_n^k:&=I_n\otimes\mathcal{L}_A^{(k)},\\
Z:&=[\,\bi{z}_1\,\bi{z}_2\,\cdots\,\bi{z}_{n-1}\,]\\
\bi{z}_i:&=\bi{e}_{ik}-\bi{e}_{ik+\bar{\kappa}}.
\end{align*}
Recall that $\bi{e}_{i}$'s are the standard basis vectors and $\otimes$ the Kronecker product defined in Section~2. Let $(\lambda_i,\bi{v}_i)$, $(\bar{\lambda}_i,\bar{\bi{v}}_i)$ and $(\tilde{\lambda}_i,\tilde{\bi{v}}_i)$ be the eigenpairs of $\mathcal{L}_A^{(k)}$, $\mathcal{L}_n^k$ and $\mathscr{L}_n^k$ respectively where $\lambda_i\leq\lambda_j$, $\bar{\lambda}_i\leq\bar{\lambda}_j$ and $\tilde{\lambda}_i\leq\tilde{\lambda}_j$ for $i<j$ and $\{\bi{v}_1,\bi{v}_2,\cdots,\bi{v}_k\}$ is an orthonormal vector set. Clearly,
$\bar{\lambda}_i=\lambda_{\lceil\frac{i}{n}\rceil}$ and we can write the modal matrix $\bar{V}$ that diagonalizes $\mathcal{L}_n^k$ as
\begin{equation}\label{def:BV}
\begin{split}
\bar{V}:&=[\,\bar{\bi{v}}_1\,\bar{\bi{v}}_2\,\cdots\,\bar{\bi{v}}_{kn}\,]\\
&=[\,I_n\otimes\bi{v}_1\quad I_n\otimes\bi{v}_2\quad \cdots\quad I_n\otimes\bi{v}_k\,].
\end{split}
\end{equation}
Before further spectral analysis of $\mathscr{L}^k_n$, we recall a classical result due to Weyl in the following.
\begin{thm}\cite[p.239]{horn13} \label{thm:weyl}
Suppose $M_1,M_2$ are Hermitian matrices of order $k$. If $\ell^{(1)}_i$, $\ell^{(2)}_i$, and $\ell^{(3)}_i$ are the $i$th smallest eigenvalues of $M_1,M_2$ and $M_1+M_2$, respectively, then, for each $i\in\{1,\cdots,k\}$,
\begin{equation} \label{eq:3-11}
	\ell^{(3)}_i \le \ell^{(1)}_{i+j} + \ell^{(2)}_{k-j},
	\quad \forall j \in \{0,1,\cdots,k-i\},
	\end{equation}
and
	\begin{equation} \label{eq:3-13}
	\ell^{(1)}_{i-j+1} + \ell^{(2)}_j \le \ell^{(3)}_i, \quad
	\forall j\in \{ 1,\cdots,i \}.
	\end{equation}	
\end{thm}
Weyl's result leads naturally to the following interlacing theorem.
\begin{cor} \label{cor:interlacing}
If in Theorem \ref{thm:weyl}, $M_2$ is in the form of $\mathbf{z}\mathbf{z}^T$ where $\mathbf{z}$ is a nonzero column vector of size $k$, then the following inequalities on the eigenvalue interlacing hold; i.e.,
	\begin{equation*}
	\begin{split}
	& \ell^{(1)}_1 \le \ell^{(3)}_1 \le \ell^{(1)}_2 \le \ell^{(3)}_2 \le \cdots \\
	& \qquad \le \ell^{(1)}_{k-1} \le \ell^{(3)}_{k-1} \le \ell^{(1)}_k \le \ell^{(3)}_k.
	\end{split}
	\end{equation*}
\end{cor}
Using Corollary~\ref{cor:interlacing} we have for
\begin{equation}\label{def:mclI}
i\in\mathcal{I}:=\{1,n+1,2n+1,\cdots,(k-1)n+1\},
\end{equation}
$\tilde{\lambda}_i=\bar{\lambda}_i$ and
\begin{equation}\label{eq:evLk}
\tilde{\bi{v}}_i^T=\left\{\begin{array}{ll}\left[\,\bi{v}^T_{\lceil\frac{i}{n}\rceil}
\,\bi{v}^T_{\lceil\frac{i}{n}\rceil}
\,\cdots\,\bi{v}^T_{\lceil\frac{i}{n}\rceil}\,\right]&\mbox{if $i=1$ or $(k-1)n+1$}\\
\left[\,\bi{v}^T_{\lceil\frac{i}{n}\rceil}\,t\bi{v}^T_{\lceil\frac{i}{n}\rceil}
\,\cdots\,t^{n-1}\bi{v}^T_{\lceil\frac{i}{n}\rceil}\,\right]&\mbox{if $i=n+1$}\\
\left[\,\bi{v}^T_{\lceil\frac{i}{n}\rceil}\,\textbf{0}^T
\,\cdots\,\textbf{0}^T\,\right]&\mbox{o.w.}
\end{array}\right.
\end{equation}
where $t=-(k-2)$. In the following we analyze the eigenpairs of $(\tilde{\lambda}_i,\tilde{\bi{v}}_i)$ for $i\notin\mathcal{I}$ in~(\ref{def:mclI}). Let $R(A,\bi{x}):=\bi{x}^TM\bi{x}/\bi{x}^T\bi{x}$ be the Rayleigh-Ritz quotient for a Hermitian matrix $M$ and nonzero column vector $\bi{x}$, and $\ell_i$ the $i$th smallest eigenvalue of $M$ with order $k$. The min-max theorem has that
\begin{align*}
\ell_i&=\min_{\dim{U} \atop=i}\max_{\bi{x}\in U}R(M,\bi{x})\\
&=\max_{\dim{U} \atop =k-i+1}\min_{\bi{x}\in U}R(M,\bi{x}).
\end{align*}
Thus for any $(k-i+1)$-dimensional $U$ we have
\begin{equation*}
\ell_i\ge\min_{\bi{x}\in U}R(M,\bi{x}),
\end{equation*}
and for any $i$-dimensional space $U$ we have
\begin{equation*}
\ell_i\le\max_{\bi{x}\in U}R(M,\bi{x} ).
\end{equation*}

\begin{lem}\label{lem:lb} For $j\in\{0,1,\cdots,k-1\}$,
\begin{equation}\label{ineq:lb}
\lambda_{j+1}=\bar{\lambda}_{nj+1}=\tilde{\lambda}_{nj+1}=\bar{\lambda}_{nj+2}<\tilde{\lambda}_{nj+2}.
\end{equation}
\end{lem}
\begin{pf} To see this, recall that for any $(nk-nj-2+1)$-dimensional $U$ we have
\begin{equation}\label{ineq:lbb}
\tilde{\lambda}_{nj+2}\ge\min_{\bi{x}\in U}R\left(\mathscr{L}_n^k,\bi{x}\right).
\end{equation}
Consider the unit vector $\bi{x}=\bar{V}\bi{y}$ where $\bar{V}$ is defined in~(\ref{def:BV}) and $\bi{y}^T:=\left[\,\textbf{o}^T\,y_{nj+2}\,y_{nj+3}\,\cdots\, y_{nk}\,\right]\in\mathbb{R}^{nk}$. Note that
\begin{align*}
R\left(\mathscr{L}_n^k,\bi{x}\right)&=\bi{y}^T\bar{V}^T\left(\bar{V}\bar{D}\bar{V}^T+ZZ^T\right)\bar{V}\bi{y}\\
&=\sum_{i=nj+2}^{nk}\bar{\lambda}_iy_i^2+\bi{y}^T\bar{V}^TZZ^T\bar{V}\bi{y}.
\end{align*}
If $y_{n(j+1)+1}=y_{n(j+1)+2}=\cdots=y_{nk}=0$, then
\begin{equation*}
R\left(\mathscr{L}_n^k,\bi{x}\right)=\bar{\lambda}_{nj+2}+
\bar{\bi{y}}^T\bar{\bar{V}}^TZZ^T\bar{\bar{V}}\bar{\bi{y}}
\end{equation*}
where
\begin{align*}
\bar{\bi{y}}:&=[\,y_{nj+2}\,y_{nj+3}\,\cdots\,y_{n(j+1)}\,]^T,\\ \bar{\bar{V}}:&=[\,\bar{\bi{v}}_{nj+2}\,\bar{\bi{v}}_{nj+3}\cdots\,\bar{\bi{v}}_{n(j+1)}\,].
\end{align*}

Observe that
$Z^T\bar{\bar{V}}$ has independent columns for each $j\in\{0,1,\cdots,k-1\}$. $\bar{\bar{V}}^TZZ^T\bar{\bar{V}}$ is positive definite and thus $R\left(\mathscr{L}_n^k,\bi{x}\right)>\bar{\lambda}_{nj+2}$. Otherwise, we have
\begin{align*}
&R\left(\mathscr{L}_n^k,\bi{x}\right)-\bar{\lambda}_{nj+2}\\
=&\sum_{i=n(j+1)+1}^{nk}\left(\bar{\lambda}_i-\bar{\lambda}_{nj+2}\right) y_j^2+\bi{y}^T\bar{V}^TZZ^T\bar{V}\bi{y}\\
\ge&\left(\bar{\lambda}_{n(j+1)+1}-\bar{\lambda}_{nj+2}\right)\sum_{i=n(j+1)+1}^{nk}y_i^2>0.
\end{align*}
By ~(\ref{ineq:lbb}), we conclude that the inequality~(\ref{ineq:lb}) holds.
\end{pf}
\begin{lem}\label{lem:ub} For $j\in\{1,2,\cdots,k-1\}$,
\begin{equation}
\tilde{\lambda}_{nj}<\tilde{\lambda}_{nj+1}=\bar{\lambda}_{nj+1}=\lambda_{j+1}.
\end{equation}
\end{lem}
\begin{pf}
The min-max theorem implies that for any $nj$-dimensional space $U$ we have
\begin{equation}
\tilde{\lambda}_{nj}\le\max_{\bi{x}\in U}R\left(\mathscr{L}_n^k,\bi{x}\right).
\end{equation}
In the following we consider a special space with dimension $nj$.
Define for $i\in\{1,2,\cdots k-1\}$,
\begin{equation}
k_i:=\left\{
\begin{array}{ll}
\sqrt{\frac{(k-2i)(k-2i+1)}{k}} &\mbox{if $i\le\lceil\frac{k}{2}\rceil-1$}\\
\sqrt{\frac{(2i-k+1)(2i-k+2)}{k}} &\mbox{o.w.}\end{array}\right.,
\end{equation}
and
\begin{equation}
p_i:=\left\{
\begin{array}{cl}
-(k-2) &\mbox{if $i=1$}\\
1 &\mbox{if $i=k-1$}\\
0&\mbox{o.w.}
\end{array}\right..
\end{equation}
Let $\mathcal{I}$ and $A_i$, $i\in\{1,2,\cdots,n\}$, be the square matrices of order $n-1$ where
\begin{equation*}
\mathcal{I}:=\left[\begin{array}{ccc}1 &\cdots & 1\\
\vdots  &\iddots &\\ 1 &&\end{array}\right]
\end{equation*}
and
\begin{equation*}
A_i:=\left[\begin{array}{ccccc}
&&&&k_i\\
&&&k_i&k_ip_i\\
&&\iddots&k_ip_i&k_ip_i^2\\
&k_i&\iddots&\vdots&\vdots\\
k_i&k_ip_i&\cdots&k_ip_i^{n-3}&k_ip_i^{n-2}
\end{array}\right].
\end{equation*}
Also, let $\bi{p}_i:=\left[\,p_i\,p_i^2\,\cdots\,p_i^{n-1}\,\right]^T$ and define the square matrices
$\mathscr{I},\tilde{\mathscr{I}},\mathcal{A}_i, P_i$ of order $n$ where
\begin{equation*}
\mathscr{I}:=\left[\begin{array}{ll}\textbf{1} &\mathcal{I}\\
1  &\textbf{0}^T\end{array}\right]\,,\quad
\tilde{\mathscr{I}}:=\left[\begin{array}{ll}\textbf{0} &\mathcal{I}\\
0  &\textbf{0}^T\end{array}\right]\,,
\end{equation*}
and
\begin{equation*}
\mathcal{A}_i:=\left[\begin{array}{ll}0 &\textbf{0}^T\\
\textbf{0}  &A_i\end{array}\right]\,,\quad
P_i:=\left[\begin{array}{ll}1 &\textbf{0}^T\\
\bi{p}_i  &\textbf{0}_{n-1}\end{array}\right].
\end{equation*}
Consider the $n(j+1)\times nj$ matrix
\begin{equation*}
W_j:=\left[\begin{array}{ccccc}\mathscr{I} &\tilde{\mathscr{I}} &\cdots &\cdots &\tilde{\mathscr{I}}\\
\mathcal{A}_1  &P_1 &&& \\
&\mathcal{A}_2 &P_2 && \\
 & &\ddots &\ddots &\\
 && &\mathcal{A}_{j-1} &P_{j-1}\\
 &&&&\mathcal{A}_j
\end{array}\right]
\end{equation*}
where $j\in\{1,2,\cdots,k-1\}$. Let the vector $\bi{x}:=\bar{V}\bar{W}_j\bi{y}$ where $\bar{V}$ is defined in~(\ref{def:BV}),
\begin{equation}
\bar{W}_j:=\left[\begin{array}{c}W_j\\\textbf{0}\end{array}\right]\in\mathbb{R}^{nk\times nj},
\end{equation}
and $\bi{y}^T=[\,y_1\,y_2\,\cdots\,y_{nj}\,]$. Observe that
$\bar{V}$ in~(\ref{def:BV}) can be written as
\begin{equation*}
\bar{V}=\left[\begin{array}{cccc}\bar{\bi{v}}_1\bi{e}_1^T&\bar{\bi{v}}_2\bi{e}_1^T
&\cdots&\bar{\bi{v}}_k\bi{e}_1^T\\
\bar{\bi{v}}_1\bi{e}_2^T&\bar{\bi{v}}_2\bi{e}_2^T
&\cdots&\bar{\bi{v}}_k\bi{e}_2^T\\
\vdots &\vdots &\ddots &\vdots\\
\bar{\bi{v}}_1\bi{e}_n^T&\bar{\bi{v}}_2\bi{e}_n^T
&\cdots&\bar{\bi{v}}_k\bi{e}_n^T\end{array}\right],
\end{equation*}
which implies $Z^T\bar{V}\bar{W}_j=\textbf{0}$ for $j\in\{1,2,\cdots,k-1\}$.
We thus have
\begin{align}
&\bi{x}^T\bi{x}-\frac{1}{\tilde{\lambda}_{nj+1}}\bi{x}^T\mathscr{L}_n^k\bi{x}\nonumber\\
=\;&\bi{y}^T\bar{W}_j\bar{V}^T\left(I_{nk}-\frac{\bar{V}\bar{D}\bar{V}^T+ZZ^T}
{\tilde{\lambda}_{nj+1}}\right)\bar{V}\bar{W}_j\bi{y}\nonumber\\
=\;&\bi{y}^T\left(W_j^TW_j-\frac{\bar{W}_j^T\bar{D}\bar{W}_j}{\tilde{\lambda}_{nj+1}}\right)\bi{y}\nonumber\\
=\;&\bi{y}^TW_j^T\left(I_{n(j+1)}-\frac{\bar{D}_{n(j+1)}}{\tilde{\lambda}_{nj+1}}\right)W_j\bi{y}\label{eq:intRR}
\end{align}
\begin{figure*}[!b]
\begin{equation}\label{eq:ltilde}
{\footnotesize
\tilde{\mathcal{L}}_A^{(k)}(\lambda)=-\left[\begin{array}{ccccccccc}
-k+\lambda	&k-1-\lambda	 &	&	&		&	&	&	&1\\
 	   &-k+1+\lambda &k-2-\lambda	 	& 	& 	& 	& 	&1	&  \\
 	   & 	&-k+2+\lambda	 &\ddots	& 	 	& 	&1	 & 	&  \\
 	 	& 	&  &\ddots &\ddots	&\iddots 	& 	& 	& \\
 	 	& 	& 	& 	&\ddots	&\ddots 	& 	& 	&\\
 	 	& 	& 	&\iddots	&	 &\ddots	&3-\lambda 	& 	&\\
 	 	& 	&1	& 	& 	& 	&-3+\lambda	&2-\lambda 	& \\
 	    1&1	& 	& 	& 	& 	& 	&-2+\lambda	& \\
      1& 	& 	& 	& 	& 	& 	&  &-1+\lambda\\
\end{array}
\right]_{k\times k}.
}
\end{equation}
\end{figure*}

where $\bar{D}_{n(j+1)}$, the upper left submatrix of $\bar{D}$, is a square matrix of order $n(j+1)$.
Note that  $I_{n(j+1)}-\frac{\bar{D}_{n(j+1)}}{\tilde{\lambda}_{nj+1}}$ is diagonal with its first $nj$ diagonal entries positive and other entries zeros. Furthermore, the special structure of $W_j$ implies that the first $nj$ rows of $W_j$ is independent, which suggests that $\left(I_{n(j+1)}-\frac{\bar{D}_{n(j+1)}}{\tilde{\lambda}_{nj+1}}\right)^{\frac{1}{2}}W_j$ has
independent columns and $(\ref{eq:intRR})$ is positive. We conclude for this particular $U$,
\begin{equation*}
\tilde{\lambda}_{nj}\leq\max_{\bi{x}\in U}\frac{\bi{x}^T\mathscr{L}_n^k\bi{x}}{\bi{x}^T\bi{x}}<\tilde{\lambda}_{nj+1}.
\end{equation*}
\end{pf}

\begin{thm}\label{thm:uni} The eigenvalues of $\mathscr{L}_n^k$ are distinct.
 \end{thm}
\begin{pf} Lemma~\ref{lem:lb} and~\ref{lem:ub} show that any number in~$\Lambda_k$ in~(\ref{def:Lda}) is a distinct eigenvalue of $\mathscr{L}_n^k$. It remains to show that the remaining eigenvalues of $\mathscr{L}_n^k$ are also distinct. Note that
\begin{align*}
\left|\mathscr{L}_n^k-\lambda I_{nk}\right|&=\left|I_n\otimes\left(\mathcal{L}_A^{(k)}-\lambda I\right)+ZZ^T\right|\\
&=\left|\left(I_n\otimes\mathcal{R}_k\right)\left(I_n\otimes\left(\mathcal{L}_A^{(k)}-\lambda I\right)+ZZ^T\right)\right|
\end{align*}
where $\mathcal{R}_k$ is a square matrix of order $k$ and its $(i,j)$th entry $r_{ij}$ satisfies
\begin{equation*}
r_{ij}:=\left\{\begin{array}{rl}1 &\mbox{if $i=j$}\\
-1 &\mbox{if $(i,j)\in\{(1,2),(2,3),\cdots,(k-2,k-1)\}$}\\
0 &\mbox{o.w.}
\end{array}\right..
\end{equation*}
That is, $\mathcal{R}_kA$ performs a sequence of row operations on the square matrix $A$ of order $k$ by replacing its $i$th row with the difference of its $i$th row and $(i+1)$th row. It turns out that
\begin{align*}
\left|\mathscr{L}_n^k-\lambda I_{nk}\right|&=\left|I_n\otimes\tilde{\mathcal{L}}_A^{(k)}(\lambda)+\tilde{Z}Z^T\right|\\
\end{align*}
where $\tilde{\mathcal{L}}^{(k)}_A(\lambda)$~(see (24)) has a similar form to that in~\cite[(9)]{hsu17b} or~\cite[(33)]{hsu18}, and $\tilde{Z}:=[\,\tilde{\bi{z}}_1\,\tilde{\bi{z}}_2\,\cdots\,\tilde{\bi{z}}_{n-1}\,]$ with
\begin{equation*}
\tilde{\bi{z}}_i:=\bi{z}_i+\bi{e}_{ik+\bar{\kappa}-1}.
\end{equation*}
Observe the $\bar{\kappa}$th row of $\tilde{\mathcal{L}}_A^{(k)}$. It is a zero row except its $\bar{\kappa}$th and $(\bar{\kappa}+1)$th entries being $\bar{\kappa}+(-1)^k-\lambda$ and $-(\bar{\kappa}+(-1)^k-\lambda)$ respectively. For the eigenpair $(\lambda,\bi{v})$ where $\lambda\notin\Lambda_k$, the $\bar{\kappa}$th and $(\bar{\kappa}+1)$th entries of $\bi{v}$ are the same number, say $c$. If $k$ is even, the $(\bar{\kappa}+1)$th row of $\tilde{\mathcal{L}}^{(k)}_A(\lambda)$
implies the $(\bar{\kappa}+2)$th entries of $\bi{v}$ is $c$ as well. In fact we can derive from the special structure of $\tilde{\mathcal{L}}^{(k)}_A(\lambda)$ that the first $k$ entries of $\bi{v}$ is $[\,c_1\,c\,c\,\cdots\,c\,c_2\,]$ for some $c_1$ and $c_2$. In case $c=0$, then the $(k-1)$row of $\tilde{\mathcal{L}}^{(k)}_A(\lambda)$ implies $c_1=0$ and thus the first row of $\tilde{\mathcal{L}}^{(k)}_A(\lambda)$ implies $c_2=0$. We can continue to argue that the last row of $\tilde{\mathcal{L}}^{(k)}_A(\lambda)$ and $\tilde{Z}$ suggest that the $(k+\bar{\kappa})$th entry of $\bi{v}$ is $0$, and thus $\bar{\kappa}$th row of $\tilde{\mathcal{L}}^{(k)}_A(\lambda)$ and $\tilde{Z}$ imply that $(k+\bar{\kappa}+1)$th entry of $\bi{v}$ is $0$ as well. Similar arguments go on to reach the conclusion that $\bi{v}=0$, a contradiction. Therefore $c$ is nonzero. Combining this result with the eigenvector property in (\ref{eq:evLk}) corresponding to the eigenvalues in $\Lambda_k$, we obtain that the $\bar{\kappa}$th and $(\bar{\kappa}+1)$th entries of any eigenvector of $\mathscr{L}_n^k$ is nonzero, and thus all eigenvalues of $\mathscr{L}_n^k$ are distinct. The case of odd $k$ can be proved similarly and is skipped.
\end{pf}
\begin{cor}\label{cor:genp} If we interconnect $n$ copies of $k$-vertex antiregular graphs by connecting the terminal vertex of the $p$th graph with the degree-repeating vertex of the $(p+1)$th graph, for each $p\in\{1,2,\cdots,n-1\}$, then the resulting graph is Laplacian controllable by an input connected to the degree-repeating vertex of the first antiregular graph.
\end{cor}

The corollary above extends the Laplacian controllability of $2n$-vertex path to its generalized version with $kn$ vertices. In the sequel we show that the case of $(2n+1)$-vertex path can be extended as well.
\begin{thm}\label{thm:onemr} Suppose $A$ is a real and square matrix of order $m$ and the first entry of every eigenvector of $A$ is nonzero. Then the first entry of every eigenvector of $\mathcal{A}$ is nonzero where
\begin{equation}
\mathcal{A}:=\left[\begin{array}{cl}0 &\textbf{0}^T\\\textbf{0} &A\end{array}\right]+(\bi{e}_1-\bi{e}_2)(\bi{e}_1-\bi{e}_2)^T
\end{equation}
\end{thm}
\begin{pf}
The condition that all the first entries of eigenvectors of $A$ are nonzero implies $A$ has distinct eigenvalues, namely, the rank of $A-\lambda I$ is $m-1$ when $\lambda$ is an eigenvalue of $A$.
We first show that all columns, excluding the first one, of $A-aI$ is independent for any $a$. In case $a$ is not an eigenvalue of $A$ then $A-aI$ is nonsingular and thus the columns are independent. If $a$ is an eigenvalue, and all but the first columns of $A-aI$ are dependent, then at least two columns of $A-aI$ are linear combinations of the remaining columns since the first column of  $A-aI$ is a linear combination of other columns. It turns out that the rank of $A-aI$ is at most $m-2$, a contradiction. Now we show that all the first entries of eigenvectors of $\mathcal{A}$ are nonzero. Suppose not. Namely, some eigenvector of $\mathcal{A}$ has a zero as its first entry. The special structure of $\mathcal{A}$ implies that the second entry of that eigenvector must be zero as well. As a result, there exists a number $c$ such that all but the first column of $A-cI$ are dependent, a contradiction.
\end{pf}
\begin{figure*}[t!]
		\centering
		\includegraphics[width=\linewidth]{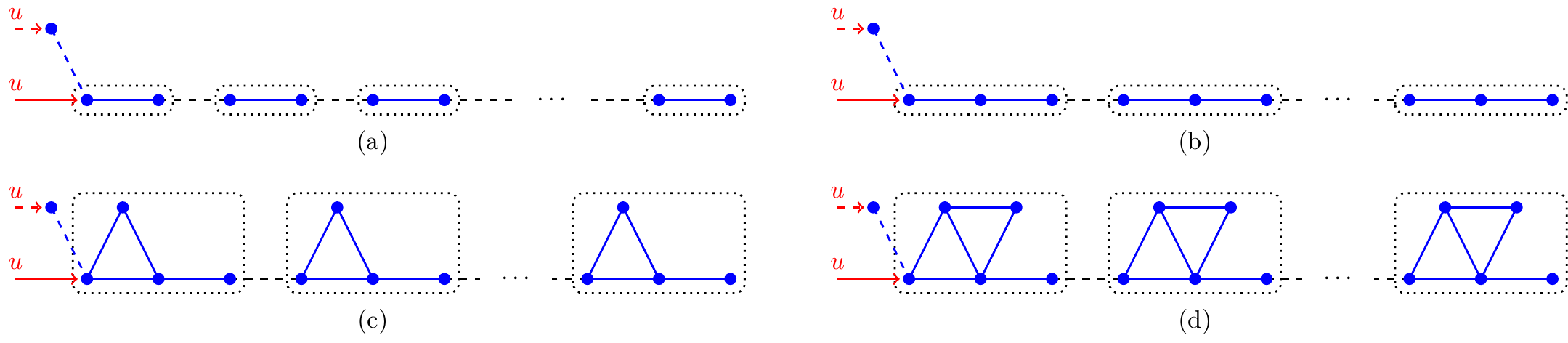}
		\caption{ Subplot (a) is a path graph whose vertex number is represented as $2n$ or $2n+1$ for some positive integer $n$. Subplots (b), (c) and (d) illustrate the generalized versions whose vertex numbers are $kn$ or $kn+1$ where $k=3,4$ and $5$ respectively. Note that if we remove the edges connecting dotted cells, the graph in any dotted cell is an antiregular graph. By Corollary~\ref{cor:genp} and~\ref{cor:xtrnode}, the control input $u$ connected to the specified node ensures the Laplacian controllability of the respective graph.}\label{fig:gen}
	\end{figure*}
Theorem~\ref{thm:onemr} allows to extend our generalization result above not only from a path with $2n$ vertices but also from $2n+1$ vertices. In the following we summarize this controllability result, whose observability counterpart has been reported in~\cite[Proposition~4.5]{par10}.

\begin{cor}\label{cor:xtrnode} If a connected simple graph is Laplacian controllable by an input connected to its node $n_c$, then the resulting graph after connecting one more node $n_a$ to node $n_c$ is Laplacian controllable by an input connected to node $n_a$.
\end{cor}
\begin{exmp} In Fig.~\ref{fig:gen} we present several examples of generalizing the Laplacian controllability of a path shown in subplot (a), where we group every two vertices with a dotted rectangle. Subplots (b), (c) and (d) increase the numbers of vertices in each dotted rectangle to $3$, $4$ and $5$ respectively. In every rectangle the vertices are connected to form an antiregular graph and these graphs are interconnected via the terminal vertex of one antiregular graph and the degree-repeating vertex of its right neighboring antiregular graph. It can be seen that the total number of vertices has the form of $kn$ or $kn+1$, where $k$'s are $2$, $3$, $4$ and $5$ in subplots (a), (b), (c) and (d) respectively. According to Corollary~\ref{cor:genp} and~\ref{cor:xtrnode}, the vertices to connect inputs to render the respective graphs Laplacian controllable are specified.
\end{exmp}

\section{Conclusions and Future Works}
We have studied in this work the signal broadcasting problem in a network. With an input that generates the signals, we would like to design a linear control system in which all state values at specific time can be specified and used to drive or operate the system. This design problem turns out to be a controllability problem of a network evolving according to the Laplacian dynamics. We have proposed a novel network topology to realize the single-input Laplacian controllability. The topology can be regarded as a generalized version of the path structure, which is one of the simplest for vertex connection. Our proposed network makes use of the single-input Laplacian controllability of an antiregular graph, and is constructed by interconnecting a finite number of these graphs in a way similar to that in constructing a path by interconnecting a finite number of two-vertex antiregular graphs. With the results we have identified a new member in the class of single-input Laplacian controllable graphs, and more importantly enriched the feasible set of single-input Laplacian controllable graphs subject to certain edge constraints. A natural issue raised by our results is that whether the property of single-input Laplacian controllability can be preserved by interconnecting antiregular graphs with different numbers of vertices? It was shown that the property is not preserved in connecting a four-vertex antiregular graph with a two-vertex path through their respective terminal vertices~\cite{hsu17b}, but it is preserved if the two-vertex path is connected to the degree-repeating vertex of the four-vertex antiregular graph~\cite{hsu18}. How to interconnect a finite number of different antiregular graphs to preserve the property is definitely worthy of a separate note. Interconnecting identical graphs other than the antiregular one, and compare their driving-force issues such as the nodal energy~\cite{zhao18} will also be a topic of interest.


\bibliographystyle{elsarticle-num}
\bibliography{thresh}
\end{document}